\documentclass{article}

\usepackage{amsmath}
\usepackage{amssymb}
\usepackage{cite}
\usepackage{graphicx}
\usepackage{xcolor}

\title{Effective Statistical Control Strategies for Complex Turbulent Dynamical Systems}
\author{Jeffrey Covington\textsuperscript{a}, Di Qi\textsuperscript{b}, Nan Chen\textsuperscript{a}}
\date{\textsuperscript{a }Department of Mathematics, University of Wisconsin-Madison, 480 Lincoln Drive, Madison, WI 53706, USA.\\ \textsuperscript{b }Department of Mathematics, Purdue University,
150 North University Street, West Lafayette, IN 47907, USA\\}

\begin{document}
    \maketitle

    \begin{abstract}
Control of complex turbulent dynamical systems involving strong nonlinearity and high degrees of internal instability is an important topic in practice. Different from traditional methods for controlling individual trajectories, controlling the statistical features of a turbulent system offers a more robust and efficient approach. Crude first-order linear response approximations were typically employed in previous works for statistical control with small initial perturbations. This paper aims to develop two new statistical control strategies for scenarios with more significant initial perturbations and stronger nonlinear responses, allowing the statistical control framework to be applied to a much wider range of problems. First, higher-order methods, incorporating the second-order terms, are developed to resolve the full control-forcing relation. The corresponding changes to recovering the forcing perturbation effectively improve the performance of the statistical control strategy. Second, a mean closure model for the mean response is developed, which is based on the explicit mean dynamics given by the underlying turbulent dynamical system. The dependence of the mean dynamics on higher-order moments is closed using linear response theory but for the response of the second-order moments to the forcing perturbation rather than the mean response directly. The performance of these methods is evaluated extensively on prototype nonlinear test models, which exhibit crucial turbulent features, including non-Gaussian statistics and regime switching with large initial perturbations. The numerical results illustrate the feasibility of different approaches due to their physical and statistical structures and provide detailed guidelines for choosing the most suitable method based on the model properties.
    \end{abstract}

    \section{Introduction}

Complex turbulent dynamical systems emerge throughout fields in science and technology, including in geophysics, engineering, neural science, and plasma physics, to name a few \cite{majda2006nonlinear, vallis2017atmospheric, majda2016introduction, nicholson1983introduction, chen2023stochastic, strogatz2018nonlinear, wilcox1988multiscale, sheard2009principles, ghil2012topics, chen2018conditional}. These turbulent systems are characterized by high-dimensional state spaces with substantial nonlinear energy transfers between scales \cite{pope2000turbulent, majda2016introduction, alexakis2018cascades}. The control of such turbulent dynamical systems has grand importance and broad applications. The goal of control is to design an optimal course of action (the control) to drive the perturbed state of interest back to a desired final target state under minimum cost within a finite time window. For example, active and passive control strategies can reduce the aerodynamic drag of vehicles such as aircraft and cargo ships \cite{sudin2014review, lee1997application, bechert1997experiments, pfeiffer2012multivariable, ahmadzadehtalatapeh2015review} and increase the efficiency of liquid and gas transport through pipelines \cite{auteri2010experimental, quadrio2011drag, kim2011physics}. Various control strategies have also been designed for applications in industrial mixing and manufacturing \cite{brodkey2012turbulence, aamo2003flow, paul2004handbook, kresta2015advances}. In addition, control of turbulent systems has significant implications for climate change mitigation involving large-scale models with high uncertainties \cite{macmartin2014dynamics, macmartin2014solar, irvine2016overview, dykema2014stratospheric}. However, controlling turbulent systems, in general, has proven a formidable challenge. The central obstacles involve the development of efficient algorithms to deal with the genuinely high dimensionality and a large number of unstable modes. Linear control theory is a well-developed field \cite{dullerud2013course, bhattacharyya2018linear, skogestad2005multivariable} and linear control strategies have been applied to turbulent systems under a variety of circumstances where the system can be linearized about a fixed-point and stabilized by the control. For example, linear control with closed-loop feedback from the system has successfully been used to delay the transition from laminar to turbulent flow \cite{bagheri2011transition, fabbiane2014adaptive, bagheri2009input}. However, linear control techniques do not scale well computationally when the dimensionality of the system becomes large \cite{kim2007linear, duriez2017machine}. Thus low-dimensional reduced-order approximations are frequently employed through model reduction and system identification \cite{willcox2002balanced, rowley2005model, brunton2016discovering, kaiser2018sparse, fukami2021sparse}. Besides linear control, there have been many other innovative approaches to the control of turbulent systems. Machine learning, for example, has been used to design control laws for turbulent systems using data \cite{duriez2017machine, bucci2019control, brunton2022data}. In addition, there have been open-loop approaches where a predetermined control is applied \cite{quadrio2011drag, van2020active}.

Statistical control offers a fundamentally different approach compared with the traditional trajectory control methods. Statistical control aims to control certain statistical features of the underlying system. These features can be considered as each statistical moment of the critical model state and are estimated by averaging an ensemble of model states. Note that even deterministic systems can fall under the statistical control framework due to the uncertainties in the initial conditions and observations. These initial uncertainties are propagated and amplified in time as a response to the strong instability considering the turbulent nature of the system. Statistical control has several unique advantages. First, there is no need to exploit exhausting procedures to resolve the full solution of each high-dimensional chaotic trajectory of the underlying turbulent dynamics. The control strategy is achieved effectively by considering only the contributions of the leading order moments. This significantly reduces the computational cost and avoids the randomness of individual trajectories in affecting the control results. Second, although individual trajectories are turbulent, the time evolution of the statistics is deterministic and is more stable. This can be seen by noticing that the Fokker-Planck equation, which is the time evolution of the probability density function (PDF) of the state variables, is always linear despite the associated underlying dynamical system being highly nonlinear and turbulent. Therefore, the statistics are more controllable. Third, statistical control can naturally account for uncertainties and incorporate stochastic reduced-order models \cite{majda2018strategies, majda2019linear}. It allows a large degree of freedom to design suitable strategies for efficient controls.

Statistical energy is a measure of the total statistical mean and variance of the system state across all scales. It is a natural scalar quantity to consider in the context of controlling turbulence \cite{farrell2014statistical, majda2015statistical, resseguier2021new}. Previous works \cite{majda2017effective, majda2019using, majda2019linear} have demonstrated that statistical responses in the key states can be successfully controlled using the statistical energy by making use of an energy-conservation principle that appears in numerous turbulent dynamical systems \cite{majda2015statistical, majda2016introduction, gugole2019numerical}. In particular, exploiting symmetry in the total statistical energy dynamics avoids the inherent nonlinear structure containing instability. Consequently, this approach eliminates the need to track and control a large number of unstable modes. In its current formulation, this statistical control strategy is run in an open-loop manner without requiring online feedback from the system, allowing for the prescribed control forcing to be determined offline for efficient computation. In addition, using the scalar-valued total statistical energy as the control object circumvents the computational issues raised by high-dimensional systems. These factors point to promising applications of statistical turbulent control considering different dynamical features of the targeting turbulent systems.

The statistical control strategy aims to control the statistical energy from a perturbed state back to the target equilibrium state by exerting an external control forcing in the underlying turbulent system. From a high-level description, the strategy consists of two consecutive steps: calculating the optimal energy control and the inversion of a nonlocal control-forcing relation. In the first step, the explicit dynamics of statistical energy are derived using the aforementioned statistical energy-conservation principle. The original control of the high-dimensional system gets reduced to a linear control problem for the scalar energy, which can then be solved using the Hamilton-Jacobi-Bellman (HJB) equations directly \cite{anderson2007optimal, bardi1997optimal}. In the second step, a nonlocal inversion problem is solved to find the deterministic external control forcing in the underlying system, which yields the optimal control of the original turbulent system. This nonlocal inversion problem uses the coupling of the optimal energy control found in the first step with both the deterministic external forcing and the response of the statistical mean of the system to the external forcing. Direct simulation of the response of the mean would be prohibitively expensive, so a crude first-order linear response approximation was employed in previous works instead \cite{majda2019linear, majda2010linear, majda2010low}. Notably, the second-order feedback term in the full control-forcing relation was truncated. This simplifies the analysis and remains consistent with the first-order approximation valid for small amplitude perturbations within the linear regime. This approach effectively controlled the statistical energy to the target equilibrium state from small initial perturbations.

This paper aims to develop new statistical control strategies for scenarios with more significant initial perturbations and stronger nonlinear responses, allowing the statistical control framework to be applied to a much wider range of problems. Note that the second-order term in the control-forcing relation, consisting of the product of the external forcing perturbation and the mean response to the forcing, is significant for large initial perturbations from the equilibrium state and thus cannot be neglected. While this term can be truncated for small perturbations, it must be included to guarantee proper performance under most large perturbations. In scenarios where the initial mean perturbation is the dominant component of the initial energy perturbation, the inclusion of the second-order term is reflected by the initial external forcing perturbation prescribed by the strategy. Even when the initial mean perturbation is small, the required strong external forcing coupled with dominant nonlinear terms to efficiently control the system usually has a correspondingly strong mean response. In such a case, the external forcing perturbation is strongly influenced by the second-order term in the control-forcing relation.

To address these difficulties, two new statistical control methods are developed in this paper. First, the higher-order methods, incorporating the second-order term, is developed to fully resolve the control-forcing relation given the mean response. The corresponding changes to recovering the forcing perturbation effectively improve the performance of the statistical control strategy in most test cases. Second, the accuracy of the mean response used in the existing statistical control methods to the external forcing also dramatically impacts the performance of the statistical control strategy. With large perturbations, although linear response theory can provide reasonable results in some special cases \cite{hairer2010simple, majda2010low, marconi2008fluctuation}, the assumptions justifying the use of linear response theory breaks down and the existing methods often lead to significant errors. In particular, the mean linear response is inadequate when the system is perturbed into a different dynamical regime than the equilibrium state. Due to these limitations, a \emph{mean closure model} for the mean response as an alternative to the mean linear response is developed in this work. The mean closure model is based on the explicit mean dynamics given by the underlying turbulent dynamical system. The dependence of the mean dynamics on higher-order moments is closed using linear response theory but for the response of the second-order moments to the forcing perturbation rather than the mean response directly. Despite still incorporating linear response theory, the introduction of explicit dynamical information from the underlying model allows the mean closure model to better reflect the properties of the perturbed regime compared to the mean linear response, which only contains information from the equilibrium statistics.

The rest of the paper is organized as follows. Section \ref{sec:review} reviews the general strategies of energy-conserving turbulent dynamical systems, including the relevant assumptions and properties needed for the statistical control strategy. 
The statistical energy control problem is formulated in Section \ref{sec:methods} along with the strategies for recovering the optimal forcing perturbation from the optimal energy control, namely combining the low-order or high-order methods with either the mean linear response or the mean closure model. In Section \ref{sec:results}, the control strategies are evaluated in detail based on two prototype models. The first model is a prototypical test model that can exhibit various behaviors and dynamical regimes. The second model is the Lorenz '96 model, which is high-dimensional and exhibits multiple dynamical regimes based on the magnitude of the external forcing. Section \ref{sec:discussion} discusses the results and provides guidance and suggestions for when the various strategies should be applied. Lastly, Section \ref{sec:conclusions} concludes the paper and offers potential future research directions.

\section{Background on Statistical Modeling} \label{sec:review}
\subsection{Statistical formulation of the turbulent systems}
Turbulent dynamical systems with quadratic energy-conserving nonlinearity can be represented in the following general canonical form \cite{majda2016introduction, majda2018strategies, majda2018model} on the state variable $\mathbf{u} \in \mathbb{R}^N$ satisfying the dynamics
    \begin{equation}
        \label{eqn:turbulent_dynamical_system} \frac{d\mathbf{u}}{dt} = (L + D)\mathbf{u} + B(\mathbf{u}, \mathbf{u}) + \mathbf{F}(t) + \boldsymbol{\sigma}(t) \dot{\mathbf{W}}(t).
    \end{equation}
Above, the linear component of the operator is decomposed into two matrices:  a skew-symmetric matrix representing linear dispersion effects, $L$, and a negative definite matrix representing dissipation effects, $D$. The quadratic nonlinearity is given by a bilinear operator, $B(\cdot, \cdot)$, which satisfies the following energy conservation law,
    \begin{equation}
        \label{eqn:energy_conservation} \mathbf{u} \cdot B(\mathbf{u}, \mathbf{u}) = 0.
    \end{equation}
Here ``$\cdot$'' denotes the Euclidean inner product. The last two terms of equation \eqref{eqn:turbulent_dynamical_system} represent the external forcing of the system, which is composed of the deterministic component of the forcing, $\mathbf{F}(t)$, and the random component, $\boldsymbol{\sigma}(t) \dot{\mathbf{W}}(t)$, where $\dot{\mathbf{W}}$ is Gaussian noise.

It is useful to decompose the state $\mathbf{u}$ into a deterministic mean state, $\bar{\mathbf{u}}(t) = \langle \mathbf{u}(t) \rangle$, and the stochastic fluctuations about each mode
    \begin{equation}
        \mathbf{u}=\bar{\mathbf{u}}+\sum_{k=1}^{N}Z_k (t)\mathbf{e}_k,
    \end{equation}
where $\langle \cdot \rangle$ denotes statistical expectation, and $\mathbf{e}_k$ is the predetermined orthonormal basis. The covariance matrix of $\mathbf{u}$ is defined as $R(t) = \langle \mathbf{Z}\mathbf{Z}^\ast\rangle$ where $\mathbf{Z} = (Z_1, \dots, Z_N)^\intercal$ and $\cdot^\ast$ denotes the conjugate transpose. Using the above mean-fluctuation decomposition of $\mathbf{u}$ and equation \eqref{eqn:turbulent_dynamical_system} the dynamics of the mean of $\mathbf{u}$ can be explicitly written as
    \begin{equation}
        \label{eqn:meanDE} \frac{d\bar{\mathbf{u}}}{dt} = (L + D)\bar{\mathbf{u}} + B(\bar{\mathbf{u}}, \bar{\mathbf{u}}) + \sum_{i,j} R_{ij}B(\mathbf{e}_i, \mathbf{e}_j) + \mathbf{F},
    \end{equation}
The operator $L_u$  incorporates the energy transfers between modes from the linear dispersion and dissipation effects
    \begin{equation}
        \{L_u\}_{ij} = \left[(L + D)\mathbf{e}_j + B(\bar{\mathbf{u}}, \mathbf{e}_j) + B(\mathbf{e}_j, \bar{\mathbf{u}})\right] \cdot \mathbf{e}_i.
    \end{equation}
	
Importantly, note that the mean dynamics given in equation \eqref{eqn:meanDE} are not closed due to the dependence on the covariance through the interactions with the nonlinearity. Further, even by including the next-order covariance dynamics from the equation. The covariance equation for $R$ can be derived as
    \begin{align}
        \label{eqn:covDE} \frac{dR}{dt} =& L_u(\bar{\mathbf{u}})R + RL_u^\ast(\bar{\mathbf{u}}) + Q_F + Q_\sigma,\\
        \label{eqn:nonlinear-energy-flux} \{Q_F\}_{ij} =& \sum_{m, n}\langle Z_m Z_n Z_j\rangle B(\mathbf{e}_m, \mathbf{e}_n) \cdot \mathbf{e}_i + \langle Z_m Z_n Z_i\rangle B(\mathbf{e}_m, \mathbf{e}_n) \cdot \mathbf{e}_j.
    \end{align}
$Q_F$ is the energy flux that accounts for the energy transfer from higher-order non-Gaussian statistics, and $Q_\sigma=\sum_{k}(\mathbf{e}_i \cdot \sigma_k)\left(\sigma_k \cdot \mathbf{e}_j\right)$ is positive definite and gives the energy transfer from the stochastic component of the external forcing. The system with \eqref{eqn:covDE} is still not closed due to the third-order moments that appear in the energy flux term $Q_F$. This means that the statistical mean response to external forcing cannot be fully resolved in this hierarchical approach, so in practice, various approximations and closures are needed \cite{majda2018strategies}.

\subsection{Statistical Energy and Response to External Forcing}

One quantity of primary interest is the total statistical energy of the system, defined as a combination of the energy in the mean and total covariance
    \begin{equation}
        E = \frac{1}{2} \bar{\mathbf{u}} \cdot \bar{\mathbf{u}} + \frac{1}{2}\operatorname{tr}(R).
    \end{equation}
Here, $\bar{\mathbf{u}}$ is the mean vector of $\mathbf{u}$, and $\operatorname{tr}(R)$ is the trace of the covariance matrix. In addition to the energy conservation principle given in equation \eqref{eqn:energy_conservation}, the following assumptions, detailed in \cite{majda2015statistical}, are needed to formulate the dynamics of the statistical energy $E$, that is,
    \begin{equation}
        B(\mathbf{e}_i, \mathbf{e}_i) \equiv 0, \quad 1 \leq i \leq N,
    \end{equation}
    and
    \begin{equation}
        \mathbf{e}_i \cdot \left[B(\mathbf{e}_j, \mathbf{e}_i) + B(\mathbf{e}_i, \mathbf{e}_j)\right] = 0 \quad \text{for any $i$, $j$.}
    \end{equation}
The above identities characterize the general symmetry in the system that the self-interactions and the closed interactions between pairs of modes vanish under the quadratic nonlinearity. To simplify the notation used in this discussion, we also assume uniform damping $D = -dI$ with $d > 0$. Under these assumptions and using equations \eqref{eqn:meanDE} and \eqref{eqn:covDE}, the total statistical energy satisfies
    \begin{equation}
        \label{eqn:energy-dynamics} \frac{dE}{dt} = -2dE + \bar{\mathbf{u}} \cdot \mathbf{F} + \frac{1}{2} \operatorname{tr}(Q_\sigma).
    \end{equation}
Statistical energy generally decays to the equilibrium state exponentially. Notably, the dynamics depends only directly on the external forcing and first-order mean state, not the covariance or higher-order moments. This allows for controlling the response to external forcing from determining the total statistical energy and solely considering the external forcing and the response of the mean to the forcing.

To determine the response of the statistical energy to perturbations of the deterministic external forcing, denote the statistical energy of the system under the equilibrium distribution as $E_\mathrm{eq}$ and denote the energy perturbation as $E'(t) = E(t) - E_\mathrm{eq}$. The equilibrium energy satisfies $dE_\mathrm{eq}/dt = 0$, so, using equation \eqref{eqn:energy-dynamics}, the equilibrium energy can be explicitly computed using only first-order mean state by
    \begin{equation}
        E_\mathrm{eq} = \frac{1}{2d} \bar{\mathbf{u}}_\mathrm{eq} \cdot \mathbf{F}_\mathrm{eq} + \frac{1}{4d} \operatorname{tr}(Q_\sigma).
    \end{equation}
Using the above equation and further denoting the deterministic forcing perturbation as $\delta \mathbf{F} = \mathbf{F} - \mathbf{F}_\mathrm{eq}$ and the corresponding mean perturbation as $\delta \bar{\mathbf{u}} = \bar{\mathbf{u}} - \bar{\mathbf{u}}_\mathrm{eq}$, the energy perturbation, $E'$, satisfies
    \begin{align*}
        \frac{dE'}{dt} = - 2dE'  + (\mathbf{F}_\mathrm{eq} \cdot \delta \bar{\mathbf{u}} + \bar{\mathbf{u}}_\mathrm{eq} \cdot \delta \mathbf{F}) + \delta \bar{\mathbf{u}} \cdot \delta \mathbf{F}.
    \end{align*}
Again, for simplicity, assuming that the stochastic component of the external forcing is not perturbed and does not contribute to the energy perturbation dynamics. It is useful to further decompose the response of the energy perturbation for each mode, i.e.
    \begin{align}
        \label{eqn:energy_response}
        \frac{dE'}{dt} &= -2dE' + \sum_{k=1}^N \left[\bar{u}_{\mathrm{eq},k} \cdot \kappa_k(t) + F_{\mathrm{eq}, k}\cdot\delta\bar{u}_k(t; \boldsymbol{\kappa}) +  \kappa_k(t) \cdot\delta\bar{u}_k(t; \boldsymbol{\kappa})\right], \\
        E'(0) &= E_0',
    \end{align}
where $\kappa_k$ is the $k$th component of $\delta \mathbf{F}$. Likewise, $F_{\mathrm{eq}, k}$, $\bar{u}_{\mathrm{eq}, k}$, and $\delta \bar{u}_{k}$ are the $k$th components of $\mathbf{F}_\mathrm{eq}$, $\bar{\mathbf{u}}_{\mathrm{eq}}$, and $\delta \bar{\mathbf{u}}$ respectively. To emphasize the central role of the forcing perturbation in each component, the mean response, $\delta \bar{u}_k$, dependence on the forcing perturbation, $\boldsymbol{\kappa} = (\kappa_1, \dots, \kappa_n)^\intercal$, is noted explicitly in the above equation.

\section{Methods on Statistical Control} \label{sec:methods}
\begin{figure}
        \begin{center}
            \includegraphics{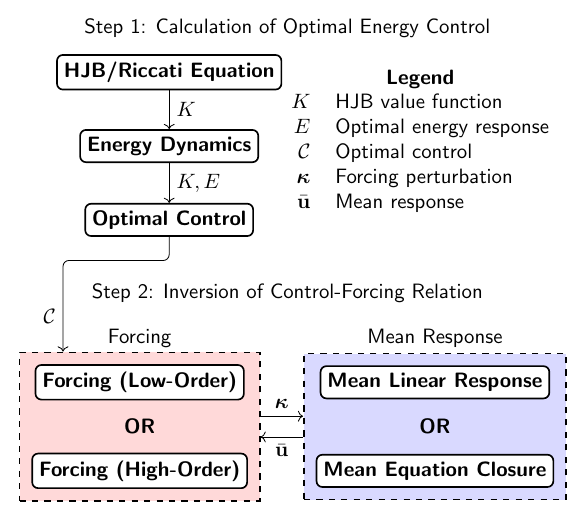}
        \end{center}
        \caption{
            Schematic diagram of the statistical control strategy. Step 1 is the calculation the optimal control, $\mathcal{C}_k$, for each mode.
            First, a Riccati equation, equation \eqref{eqn:riccati_equation}, is solved.
            This is used to calculate the optimal energy response using equation \eqref{eqn:optimal_energy}.
            The optimal control is then calculated using equation \eqref{eqn:optimal_control}.
            Step 2 consists of finding the forcing perturbation, $\kappa_k$, which yields the optimal control in each mode.
            Inverting the control-forcing relation involves solving coupled equations for the forcing and the mean response to that forcing.
            There are two choices for the forcing equations: the low order equations, equation \eqref{eqn:low_order_method}, and the high order equations, equation \eqref{eqn:high_order_method}.
            For the mean response there are two strategies: a linear response for the mean, equation \eqref{eqn:mean_linear_response}, and the mean dynamics with a closure for the higher-order moments, equation \eqref{eqn:mean_closure}.
            Choosing one strategy from each category yields four strategies total.
        }\label{fig:ctrl_diagram}
    \end{figure}

This section describes the control strategies for high-dimensional turbulent systems, including large-amplitude perturbations. The method is generally split into two consecutive steps: i) the recovery of the optimal energy control for the total statistical energy, and ii) the attribution of the forcing contribution for each detailed spectral mode. Especially, high-order accuracy is achieved by considering different ways to approximate the mean responses and high-order feedback in the control. We illustrate the general idea in the diagram in Figure~\ref{fig:ctrl_diagram}.

\subsection{Optimal Control of the Perturbed Energy}

As the first step of the statistical control strategy, the objective is to drive the total statistical energy from a perturbed state back to a target equilibrium state with a minimized cost through prescribing the deterministic external forcing perturbation in each mode, $\kappa_k$. The direct optimal control of the energy dynamics is achieved by controlling the much simpler scalar equation independent of the full dimensionality of the system. In the second step, the external forcing perturbation that yields this optimal control will be recovered by attributing the contribution to the total energy from each individual mode.

In this and future sections, the energy perturbation will be denoted as $E$ to simplify the notation. Define the energy control problem as
    \begin{equation}
        \frac{dE}{dt} = -2dE + \sum_{k=1}^N \mathcal{C}_k
    \end{equation}
where the objective is to control the energy perturbation, $E$, back to zero over the time interval $[0, T]$ using the controls $\mathcal{C}_k$, representing the total contribution to the energy perturbation from each mode. The cost functional of this control problem is proposed as
    \begin{equation}
        \label{eqn:cost}
        \mathcal{F}_\alpha[\mathcal{C}_k(\cdot)] \equiv \int_t^T \left[E^2(s) + \sum_{k=1}^M \alpha_k \mathcal{C}_k^2(s)\right] \,ds + k_T E^2(T),
    \end{equation}
where $\alpha_k$ gives the relative weights between each control and the total energy perturbation. $k_T$ is the cost coefficient for the final energy perturbation from the equilibrium state at time $T$. From the above setups, the control of total energy becomes a standard linear control problem with a quadratic cost, and so the Hamilton-Jacobi-Bellman (HJB) equation \cite{anderson2007optimal, bardi1997optimal} can be applied to find the optimal control $\mathcal{C}_k^\ast$. In addition, the total statistical energy is a scalar quantity, so the associated energy control problem is tractable even when the underlying system is high-dimensional.

Solving the HJB equations \cite{majda2019using} leads to the following Riccati equation
    \begin{align}
        \label{eqn:riccati_equation}
        \frac{dK}{dt} =& \sum_{k} \alpha_k^{-1}K^2 + 4dK - 1, \quad 0 \leq t < T \\
         K(T) =& k_T
    \end{align}
which is solved backward in time. The quantity $K$ is a factor of the value function that appears in the HJB equations. The full details of this application of the HJB equations can be found in \cite{majda2019using}. Using the solution of $K$, the optimal response of the energy perturbation, $E^\ast$, is given by the forward equation
    \begin{align}
        \label{eqn:optimal_energy}\frac{dE^\ast}{dt} =& -\left(2d + \sum_{k} \alpha_k^{-1} K\right)E^\ast, \quad 0 \leq t < T, \\
        E(0) =& E_0.
    \end{align}
Finally, the optimal control in each mode, $\mathcal{C}_k^\ast$, can be calculated as
    \begin{align}
        \label{eqn:optimal_control} \mathcal{C}_k^\ast(t) =& -\alpha_k^{-1}K(t)E^\ast(t).
    \end{align}

\subsection{Inversion of the Control-Forcing Relation}

The goal is to find the forcing perturbation in each mode $\kappa_k$ that yields the optional control discovered from the energy perturbation. For simplicity in notation, the optimal control of each mode found in the previous section will be denoted by $\mathcal{C}_k$, where the superscript ``$\ast$'' is dropped. Using equation \eqref{eqn:energy_response} yields the following relation between the energy control and the forcing perturbation in each mode:
    \begin{equation}
        \label{eqn:control_forcing_relation}
        \mathcal{C}_k(t) = \bar{u}_{\mathrm{eq},k} \cdot \kappa_k(t) + F_{\mathrm{eq}, k}\cdot\delta\bar{u}_k(t; \boldsymbol{\kappa}) +  \kappa_k(t) \cdot\delta\bar{u}_k(t; \boldsymbol{\kappa}).
    \end{equation}
The mean perturbation response, $\delta\bar{u}_k$, depends on the forcing perturbation $\boldsymbol{\kappa}$. Inverting this relation for $\kappa_k$ involves solving an ODE system which couples the optimal control, $\mathcal{C}_k$, the forcing perturbation, $\kappa_k$, and the mean response, $\delta\bar{u}_k$.

\subsubsection{The Low-Order Method}
In previous works, small forcing perturbation and mean state response are assumed. Thus, the second-order term, $\kappa_k \cdot \delta\bar{u}_k$, is omitted in equation \eqref{eqn:control_forcing_relation} so that the relation only accounts for the dominant leading-order response of the energy to the forcing perturbation \cite{majda2017effective, majda2019using}. This assumption is justified by using leading order approximations for the mean response for small initial perturbations, which introduce $O(\delta^2)$ errors in the same order as the second-order perturbation term.  In doing so, the inversion relation is linearized, leading to simplified analytical solutions. The ODE resulting from this linearized relation, referred to in this paper as the ``low-order method'', is given by
    \begin{align}
        \label{eqn:low_order_method}
        \frac{d \kappa_k}{dt}  =& \frac{1}{u_{\mathrm{eq}, k}} \left(\frac{d\mathcal{C}_k}{dt} - F_{\mathrm{eq}, k} \cdot \frac{d\bar{u}_k}{dt} \right) \\
        \kappa_k(0) =& \frac{1}{\bar{u}_{\mathrm{eq}, k}} \left(\mathcal{C}_k(0) - F_{\mathrm{eq}, k} \cdot \delta \bar{u}_k(0)\right).
    \end{align}
Here $d\mathcal{C}_k/dt$ can be explicitly calculated using equations \eqref{eqn:riccati_equation}, \eqref{eqn:optimal_energy}, and \eqref{eqn:optimal_control}
    \begin{equation}
        \label{eqn:control_derivative}
        \frac{d\mathcal{C}}{dt} = -\alpha_k^{-1} E^\ast(t) (2dK(t) - 1).
    \end{equation}
The term $\delta \bar{u}_k(0)$ denotes the mean perturbation in the initial perturbed state. Solving this ODE system involves approximating the response of the mean, $d\bar{u}_k/dt$, to the forcing perturbation. Strategies for computing the mean responses are detailed in sections \ref{subsec:mean_resp}.

\subsubsection{The High-Order Method}

The second-order term in the control forcing relation was truncated in the low-order method. However, large errors might be introduced due to this omission when the perturbation amplitude grows large. Still, the same analysis can be done by including this additional term. The ODE resulting from inverting equation \eqref{eqn:control_forcing_relation}, including all contributing terms, is given by
    \begin{align}
        \label{eqn:high_order_method}
        \frac{d \kappa_k}{dt}  =& \frac{1}{u_{\mathrm{eq}, k} + \delta \bar{u}_k(t)} \left(\frac{d\mathcal{C}_k}{dt} - \left(F_{\mathrm{eq}, k} + \kappa_k(t)\right) \cdot \frac{d\bar{u}_k}{dt} \right) \\
        \kappa_k(0) =& \frac{1}{\bar{u}_{\mathrm{eq}, k} + \delta \bar{u}_k(0)} \left(C_k(0) - F_{\mathrm{eq}, k} \cdot \delta \bar{u}_k(0)\right),
    \end{align}
where $d\mathcal{C}_k/dt$ is given by equation \eqref{eqn:control_derivative} and $\delta \bar{u}_k(0)$ is the initial mean perturbation in the perturbed state. Compared to equation \eqref{eqn:low_order_method}, equation \eqref{eqn:high_order_method} contains one more perturbation term that accounts for the higher-order contributions to the energy response. For large perturbations from the target equilibrium state, the high-order feedback term becomes necessary to accurately recover the true forcing forms. The inversion of the full control-forcing relation will be referred to in this paper as the ``high-order method'' in which the response of the energy to the forcing and mean perturbations is fully resolved. On the other hand, the extra terms in the denominator may lead to additional numerical complications, but the improved performance is generally worth the additional sophistication. We will discuss the performance in the detailed numerical tests in Section~\ref{sec:results}.

\subsection{Different Strategies to Recover Mean Responses}\label{subsec:mean_resp}



The inversion of the control-forcing relation given by equations \eqref{eqn:low_order_method} and \eqref{eqn:high_order_method} requires the response of the mean to the forcing perturbation. Unfortunately, a direct simulation for the mean responses would be prohibitively expensive considering the extremely high dimensional problem, so approximations for the mean response are developed instead. In this subsection, we develop two approaches to efficiently estimate the mean responses without directly solving the full equations. The first strategy uses linear response theory as a convenient way to recover the leading-order mean response based on the Fluctuation-Dissipation Theorem (FDT) \cite{majda2005information, leith1975climate}. The second approximation provides a more accurate high-order approximation based on solving the explicit mean dynamical equation with a high-order closure.

\subsubsection{Mean Linear Response to Forcing} \label{sec:meanlinear}
The linear mean response to the forcing perturbation $\boldsymbol{\kappa}$ is computed by
    \begin{equation}
        \label{eqn:mean_linear_response}
        \delta \bar{u}_k = \sum_{\ell = 1}^N \left[\int_0^t \mathcal{R}_{\bar{u}, k \ell} (t - s) \kappa_\ell(s) \, ds + \int_{-\infty}^0 \mathcal{R}_{\bar{u}, k \ell} (t - s) \delta F_{\mathrm{p}, \ell}\, ds\right] + O(\delta^2),
    \end{equation}
where the \emph{linear response operator} is given by
    \begin{equation}
        \label{eqn:linear_response_operator}
        \mathcal{R}_{\bar{u}, k\ell}(t) = \langle (\bar{u}_k(t) - \bar{u}_{\mathrm{eq}, k}) G_\ell[\bar{\mathbf{u}}(0)] \rangle_\mathrm{eq},
    \end{equation}
    and
    \begin{equation}
        G_\ell(\mathbf{u}) = -p_\mathrm{eq} \operatorname{div}_\mathbf{u}(\mathbf{e}_\ell p_\mathrm{eq}(\mathbf{u})),
    \end{equation}
where $p_\mathrm{eq}(\mathbf{u})$ is the equilibrium probability density of $\mathbf{u}$. The basis vector $\mathbf{e}_\ell$ corresponds with the $\ell$th mode. The forcing perturbation $\delta\mathbf{F}_{\mathrm{p}, \ell}$ is the constant external forcing perturbation corresponding to the initial perturbed state. This is contrasted with $\boldsymbol{\kappa}$, which is the time-dependent forcing perturbation based on the control that takes over at time $t = 0$. The linear response operator only depends on the PDF of the target equilibrium state. This explicit formula provides a great computation reduction for convenient calculation of the mean responses but only in the leading order.

Then, the linear response approximation is incorporated into the low-order and high-order methods given in equations \eqref{eqn:low_order_method} and \eqref{eqn:high_order_method}. The equations of the mean response are given by
    \begin{equation}
        \label{eqn:mean_linear_response_deriv}
        \frac{d \bar{u}_k}{dt} = \sum_{\ell = 1}^N \left[\mathcal{R}_{k, \ell}(0) \kappa_{\ell}(t) + \int_0^t \mathcal{R}_{\bar{u}, k \ell}'(t - s) \kappa_\ell(s) \, ds - \mathcal{R}_{k, \ell}(t) \delta F_{\mathrm{p}, \ell}\right]
    \end{equation}
    with initial condition
    \begin{equation}
        \delta \bar{u}_k(0) = \sum_{\ell = 1}^N \left[\int_{-\infty}^0 \mathcal{R}_{\bar{u}, k \ell} (t - s) \delta F_{\mathrm{p}, \ell}\, ds\right].
    \end{equation}
While the linear response includes contributions from all modes, typically, only the contribution from a few modes is relevant, allowing the rest to be truncated to further save the computation cost.

As further notice, the linear response operator can still be difficult to calculate. In this paper, we adopt the quasi-Gaussian approximation \cite{leith1975climate, majda2010high, gershgorin2010test} where the equilibrium probability density, $p_\mathrm{eq}(\mathbf{u})$, is approximated by a Gaussian distribution. In this case
    \begin{equation}
        G_\ell(\mathbf{u}) = \mathbf{e}_\ell \cdot \mathbf{R}_\mathrm{eq}^{-1} (\mathbf{u} - \bar{\mathbf{u}}_\mathrm{eq})
    \end{equation}
is used in equation \eqref{eqn:linear_response_operator}. Notice this approximation is quasi-Gaussian since the higher-order moments will still be involved due to the probability expectation in \eqref{eqn:linear_response_operator} based on the true solution rather than only the Gaussian closure.

There are many other strategies for approximating the linear response operator \cite{majda2019linear}. While the quasi-Gaussian approximation is more than sufficient for our purposes, other strategies can be investigated. The previous works on the statistical control strategy have used an exponential fit based on an information theory criterion to approximate the linear response operator
    \begin{equation}
        \mathcal{R}_{k\ell}(t) \approx \exp(-\gamma_k t).
    \end{equation}
This is justified with a quasi-Gaussian approximation and a diagonal covariance matrix where the linear response operator reduces to an autocorrelation function of $\mathbf{u}$, which frequently have exponential and oscillatory structures. While this exponential fit is not utilized in this paper, it is noted for its potential application as a useful approximation in practical settings and for the theoretical convenience of having an explicit form of the linear response operator.

\subsubsection{Mean Dynamical Equation Closure for Response to Forcing} \label{sec:meanclosure}

While directly using linear response theory for the mean response provides a useful approximation, as in section \ref{sec:meanlinear}, the linear response is proved inadequate in many cases with large perturbations. This paper develops another strategy to incorporate the mean response based directly on the mean dynamical equation explicitly given in equation \eqref{eqn:meanDE} subject to the forcing perturbation $\boldsymbol{\kappa}$. For this method, the dependence of the mean dynamics on the second-order moments through the nonlinearity is closed using a linear response for the covariance. Compared to the linear response of the mean, which only incorporates statistical information from the target equilibrium distribution, the use of the mean dynamical equation incorporates crucial additional dynamical information. This can be particularly useful when the system is perturbed into a different dynamical regime.

The mean dynamical equation in the mean closure method is given by \eqref{eqn:meanDE}
    \begin{align}
        \label{eqn:mean_closure}
        \frac{d\bar{\mathbf{u}}}{dt} =& (L + D)\bar{\mathbf{u}} + B(\bar{\mathbf{u}}, \bar{\mathbf{u}}) + \sum_{i,j} R_{ij}(\boldsymbol{\kappa})B(\mathbf{e}_i, \mathbf{e}_j) + \mathbf{F}_\mathrm{eq} + \boldsymbol{\kappa}.
    \end{align}
This equation still requires the solution of the second-order covariances $R$ inversely dependent on the form of the forcing perturbation $\boldsymbol{\kappa}$. A closure model is constructed using linear response theory for the response of the covariance, $R_{ij}$, to the external forcing
    \begin{multline}
        \label{eqn:covariance_linear_response}
        R_{ij}(t;\boldsymbol{\kappa}) = R_{\mathrm{eq}, ij} + \sum_{\ell=1}^N \bigg[ \int_0^t \mathcal{R}_{R, ij\ell}(t - s)\kappa_{\ell}(s) \, ds \\ + \int_{-\infty}^0 \mathcal{R}_{R, ij\ell}(t - s) \delta F_{\mathrm{p}, \ell} \, ds \bigg] + O(\delta^2),
    \end{multline}
where the linear response operator for the covariance is given by
    \begin{align}
        \mathcal{R}_{R, ij\ell}(t) = \langle (\bar{u}_i(t) - \bar{u}_{\mathrm{eq}, i})(\bar{u}_j(t) - \bar{u}_{\mathrm{eq}, j}) G_\ell[\bar{\mathbf{u}}(0)] \rangle_\mathrm{eq},\label{eqn:cov_resp}
    \end{align}
    and
    \begin{align*}
        G_\ell(\mathbf{u}) = -p_\mathrm{eq} \operatorname{div}_\mathbf{u}(\mathbf{e}_\ell p_\mathrm{eq}).
    \end{align*}
The higher-order closure of the mean equation enables a better characterization of the mean responses respecting its explicit nonlinear dynamics. The errors from the linear response approximation then appear in the second-order covariances rather than the first-order mean. The linear response in \eqref{eqn:cov_resp} will require the computation of lagged third moments, adding more non-Gaussian information into the approximation. As in Section~\ref{sec:meanlinear}, a quasi-Gaussian approximation is used for the linear response operator to efficiently compute the response operators.

\section{Numerical Results} \label{sec:results}

We examine the performance of the methods developed in Section~\ref{sec:methods} using detailed numerical tests. Combining the high-order or low-order methods with the mean equation closure or mean linear response gives a total of four approaches to compare. These methods concern different aspects of component approximations of the statistical control strategy, and the choice to implement each one can be made accordingly based on the specific problem.

These four strategies are evaluated on two complex nonlinear models exhibiting various dynamical and statistical behaviors. The first test model is a prototypical triad nonlinear model \cite{majda2018strategies} focusing on a generic coupling between three modes of a turbulent system. It can exhibit a wide variety of nonlinear and non-Gaussian behaviors. Two regimes of this model are considered, including a nearly Gaussian regime with nonlinear energy transfers between modes and a non-Gaussian regime exhibiting an energy cascade representing the transition to turbulence. The second test model is the classic Lorenz '96 model \cite{lorenz1996predictability} with 40 dimensions which shows multiple dynamical regimes depending on the magnitude of the external forcing. Large perturbations inducing regime switching will be the primary consideration. These test models will illustrate the differences between the statistical control strategies.

The experimental setup is as follows. First, the external forcing perturbation is calculated offline using different statistical control strategies. The optimal energy control is calculated using equations \eqref{eqn:riccati_equation}, \eqref{eqn:optimal_energy}, and \eqref{eqn:optimal_control}. The forcing perturbation which yields this optimal control is found using either the high-order method given in equation \eqref{eqn:high_order_method} or the low-order method given in equation \eqref{eqn:low_order_method}. The mean response to the forcing perturbation is approximated by the mean equation closure model described in equations \eqref{eqn:mean_closure} and \eqref{eqn:covariance_linear_response} or by the mean linear response in equations \eqref{eqn:mean_linear_response} and \eqref{eqn:mean_linear_response_deriv}. Second, the external forcing perturbation is applied to a Monte Carlo simulation of the underlying dynamical system. An initial ensemble of model trajectories of size $M=1\times10^4$ is drawn from the initial distribution. The perturbed initial state is created by a deterministic forcing perturbation. The deterministic component of the external forcing $\mathbf{F}_\mathrm{eq}$ is perturbed by a constant forcing amplitude $\delta \mathbf{F}_\mathrm{p}$ to drive the ensemble is forced into a perturbed state away from the original statistical equilibrium.

At the time $t=0$, the control strategy takes over. Thus the previous forcing perturbation is replaced by the statistical control forcing, $\boldsymbol{\kappa}$. The response of the statistical energy to the forcing is calculated from the ensemble and is tracked as the system is controlled back to the equilibrium state. The energy responses from each strategy are compared to the theoretically optimal energy response in equation \eqref{eqn:optimal_energy}. The uncontrolled case where no control forcing perturbation is applied, in which the energy naturally decays back to the original equilibrium state, is also used as a point of comparison. Other quantities, such as the mean response, variance response, and empirical control, are also compared for tracking the performance.

    \subsection{A Prototype Nonlinear Triad Model}\label{sec:triad}


    \begin{figure}
        \begin{center}
            \includegraphics{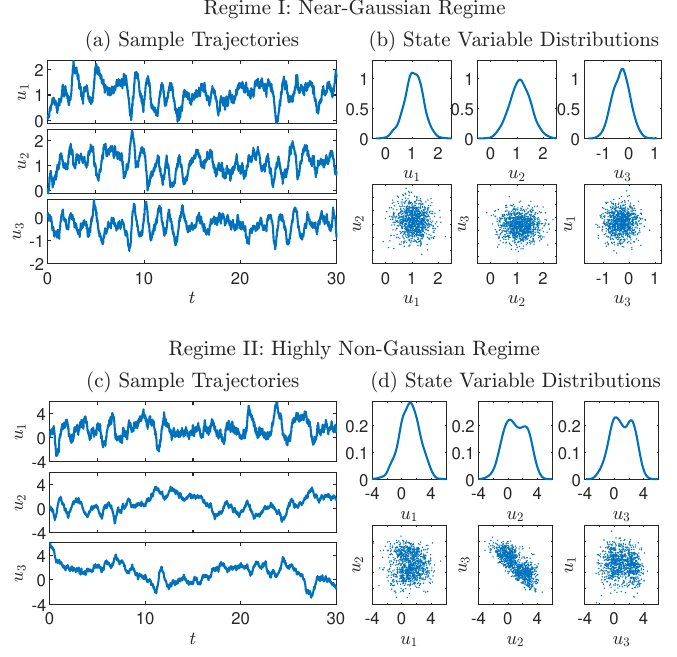}
        \end{center}
        \caption{\label{fig:triad_regimes}
            The dynamics and equilibrium distributions of the prototype triad model under two different regimes.
            Panels (a) and (c) show sample trajectories of each regime of the model.
            Panels (b) and (d) show the equilibrium marginal distributions of the state variables as well as their pairwise joint distributions in each regime.
            The nonlinearity in the model produces non-Gaussian distributions in each regime.
            In particular Regime II exhibits intermittency and highly non-Gaussian statistics.
        }
    \end{figure}

The first test model is a prototypical 3-dimensional model equipped with a quadratic energy-conserving nonlinearity, which is referred to in this paper as the \emph{triad model}. The triad model represents a generic nonlinear coupling between three variables universal in turbulent flows \cite{majda2018strategies, majda2012physics, chen2018rigorous}. Such a triad interacting structure would emerge as the bare truncation of three identified modes from a high-dimensional turbulent model. The triad model can also generate a wide variety of nonlinear and non-Gaussian behaviors due to the dominant role of the nonlinear coupling term. This sets a desirable first test model to evaluate the skills of the different proposed approaches considering the high-order contributions in the mean state and the energy equation. Despite the nonlinearity, the triad model is analytically tractable in terms of the equilibrium statistics when the linear parts have certain special structures \cite{majda2012physics}, making it an appropriate test model used to have an in-depth study of the different features of the four proposed strategies.

The state variables of the triad model are represented by $\mathbf{u} = (u_1, u_2, u_3)^\intercal$ with governing differential equations:
    \begin{align}
        \frac{du_1}{dt} =& L_2 u_3  - L_3 u_2 -d_1 u_1 + B_1 u_2 u_3 + F_1 + \sigma_1 \dot{W}_1, \\
        \frac{du_2}{dt} =& L_3 u_1  - L_1 u_3 - d_2 u_2 + B_2 u_3 u_1 + F_2 + \sigma_2 \dot{W}_2, \\
        \frac{du_3}{dt} =& L_1 u_2 -L_2 u_1  - d_3 u_3 + B_3 u_1 u_2 + F_3 + \sigma_3 \dot{W}_3.
    \end{align}
In addition, the quadratic coupling coefficients satisfy
    \begin{equation}
        B_1 + B_2 + B_3 = 0,
    \end{equation}
which ensures the general energy-conserving property \eqref{eqn:energy_conservation} of the turbulent dynamical system framework. Figure \ref{fig:triad_regimes} shows two typical regimes of the triad model used to evaluate the strategies: one has near-Gaussian statistics, while the other is highly non-Gaussian.

\subsubsection{Control on a Near-Gaussian Regime}
    \begin{figure}
        \begin{center}
            \includegraphics{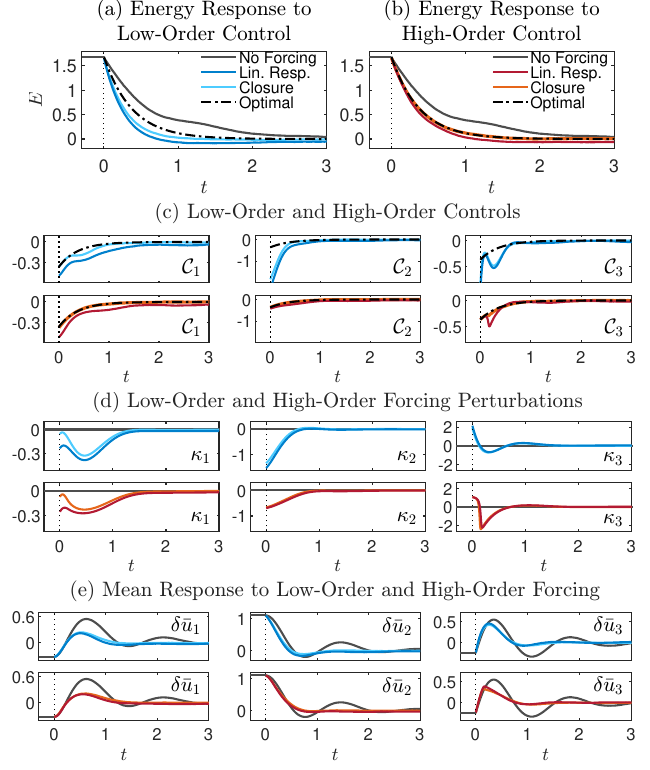}
        \end{center}
        \caption{\label{fig:ep_regime}
            The control of the prototype triad model from the perturbed state back to the equilibrium state in the near-Gaussian regime.
            Panels (a) and (b) show the energy response to the forcing, including the response from no control, the linear response and equation closure strategies, as well as the optimal response.
            Panel (a) shows the energy response for the low-order strategies: using a mean linear response and using a mean equation closure model.
            Panel (b) shows the same with the high-order strategies.
            Panel (c) compares the controls realized by the various strategies to the optimal control.
            Panel (d) shows the forcing perturbations prescribed by each strategy.
            Panel (e) shows the responses of the mean under each strategy.
        }
    \end{figure}

The first test regime for the triad model is a near-Gaussian regime which nonetheless contains strong nonlinear energy transfers between modes to reach the equipartition of energy. Sample trajectories and equilibrium distributions for this regime are pictured in Regime I in Figure \ref{fig:triad_regimes}. Nearly Gaussian and weakly non-Gaussian features are common in practice, such as in fully turbulent flow with strong mixing. The damping coefficients for this regime are  $d_1 = d_2 = d_3 = 1$. The linear dispersion coefficients are $L_1 = 3$, $L_2 = 2$, and $L_3 = -1$. The nonlinear quadratic coupling coefficients are $B_1 = 1$, $B_2 = -0.6$, and $B_3 = -0.4$. The deterministic external forcing for the equilibrium state is given by $F_1 = F_2 = 1$ and $F_3 = -1$ while the stochastic external forcing coefficients are given by $\sigma_1 = \sigma_2 = \sigma_3 = 0.5$. To perturb the model state, $F_3$ is perturbed by $\delta F_{\mathrm{p},3} = -4$ until time $t = 0$. While only $F_3$ is perturbed initially, all modes are used to control the system back to the equilibrium state.

The control of the system from the perturbed state back to the equilibrium state under different strategies is shown in Figure \ref{fig:ep_regime}. All the methods show faster convergence than the no-control scenario to efficiently return the unperturbed equilibrium state. Notably, the high-order method achieves the most accurate near-optimal performance, particularly the high-order method with a mean equation closure. In contrast, the low-order methods overshoot the response incurring relatively higher costs. This is the first confirmation of the crucial role of the higher-order correction in the energy equation when nonlinearity is dominant, as in the triad system. There is a stark difference in the control forcing perturbations between the high-order and low-order methods which can be seen by $\kappa_2$ and $\kappa_3$ in panel (d) of Figure \ref{fig:ep_regime}. The initial forcing perturbation for $\kappa_2$ differs entirely from the low-order and high-order methods. This is because the corresponding initial mean perturbation, $\delta \bar{u}_2$, pictured in panel (e), is relatively large; thus, the initial contribution of the second-order $\kappa_2 \cdot \delta \bar{u}_2$ term to the control-forcing relation is quite significant. Indeed, the low-order method overcompensates for lacking this high-order term with a large initial forcing perturbation. In contrast, the initial forcing perturbation for the high-order method is significantly smaller. For the third mode, $\kappa_3$, the initial forcing perturbation is comparable between the high-order and low-order methods due to the initial mean perturbation $\delta \bar{u}_3(0)$ being relatively small. However, the high-order method produces a stronger forcing perturbation in $\kappa_3$ shortly after the initial time. This is explained by the observation that in all methods, the mean perturbation response, $\delta \bar{u}_3$, quickly takes on a large value. So the second-order $\kappa_3 \cdot \delta \bar{u}_3$ term in the control-forcing relation is non-negligible. Because of this, only the high-order method can account for the contribution of the second-order term to the energy response. One can also see a difference between the initial forcing perturbation for $\kappa_1$ between the mean linear response and mean equation closure methods. This is an example where the mean linear response does not accurately produce the initial mean perturbation, while the mean closure model directly incorporates the initial mean perturbation.

	\subsubsection{Control on a highly non-Gaussian Regime}
	
	    \begin{figure}
        \begin{center}
        \includegraphics[scale=1.2]{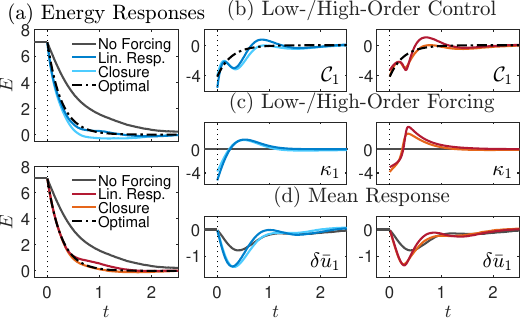}
        \end{center}

        \caption{
            \label{fig:nongaussian_regime}
            Example of controlling a highly non-Gaussian regime in the prototypical triad model.
            Panel (a) shows the response of the energy to the forcing perturbation for all strategies.
            Panel (b) shows the control for each strategy for the $u_1$ mode.
            Panel (c) shows the forcing perturbation in the $u_1$ mode.
            Panel (d) shows the mean response in the $u_1$ mode.
        }
    \end{figure}

The second test regime, shown in Regime II of Figure \ref{fig:triad_regimes}, has highly non-Gaussian statistics and intermittency. It features an energy cascade from $u_1$ to $u_2$ and $u_3$ reminiscent of the transition to turbulence. The damping coefficients are $d_1 = d_2 = d_3 = 1$, the quadratic nonlinear coupling coefficient are $B_1 = 2$ and $B_2 = B_3 = -1$, and the linear dispersion coefficients are $L_1 = 0.03$, $L_2 = 0.02$, and $L_3 = -0.01$. The unperturbed deterministic external forcing is given by $F_1 = F_2 = F_3 = 2$ and the stochastic external forcing coefficients are $\sigma_1 = 2$ and $\sigma_2 = \sigma_3 = 1$. The perturbed state is achieved through constant deterministic forcing perturbations $\delta F_{\mathrm{p}, 1} = \delta F_{\mathrm{p}, 2} = \delta F_{\mathrm{p}, 3} = 2$ until time $t = 0$. Figure \ref{fig:nongaussian_regime} shows the results of applying the control strategies to Regime II of the triad model. We focus on the performance of the dominant mode $u_1$. The other two modes $u_2 u_3$ have qualitatively similar performance and are omitted for a cleaner representation. Similar to the near-Gaussian regime, there is a strong forcing perturbation in $\kappa_1$ for the high-order method after the initial time. The initial forcing using different methods is comparable due to the relatively small initial mean perturbation $\delta \bar{u}_1(0)$ but the mean response shortly after the initial time requires the high-order method to capture the subsequent response by the energy. In addition, this example illustrates how the mean equation closure model can produce more accurate forcing perturbations under a strongly nonlinear non-Gaussian regime even when the initial perturbation is similar. As expected, stronger non-Gaussianity requires a more accurate calibration of the mean responses taking into account the higher-order statistics. The high-order equation closure gains a more accurate estimation of the mean state, thus leading to the most accurate result. For low-order methods, lacking the higher-order correction term often leads to larger errors. We suspect that the agreement in the low-order linear response approach comes as an accidental cancellation of errors.

\subsection{A High-Dimensional Model with Multiple Regimes} \label{sec:lorenz}

    \begin{figure}
        \begin{center}
            \includegraphics{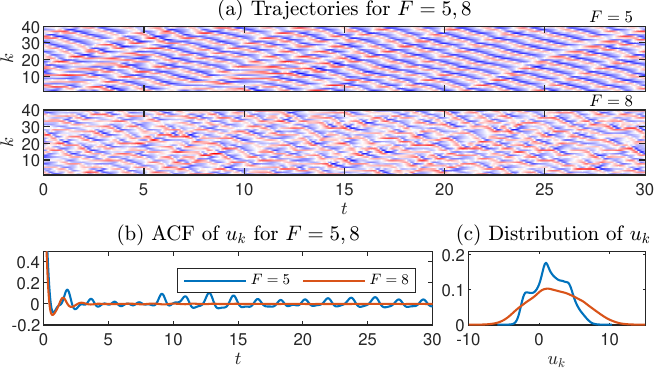}
        \end{center}
        \caption{\label{fig:lorenz_dynamics}
            Sample trajectories and distributions of the 40-dimensional Lorenz 96 model for both the $F = 5$ (weakly chaotic; highly non-Gaussian) and $F = 8$ (strongly chaotic; nearly Gaussian) regimes.
            Panel (a) shows sample trajectories for one sample of each regime in the form of the Hovmoller diagram.
            Panel (b) shows the autocorrelation functions (ACFs) for each regime.
            Note that the ACF for the $F = 5$ regime exhibits long-term oscillatory behavior while the ACF of the $F = 8$ regime decays very fast.
            Panel (c) shows the equilibrium distribution for each regime. The $F = 5$ regime is highly non-Gaussian while the $F=8$ regime is nearly Gaussian.
        }

    \end{figure}

The Lorenz '96 model is a standard test model which mimics geophysical waves and exhibits phenomena such as mid-latitude baroclinic instability \cite{lorenz1996predictability}.    The model is defined in a 40-dimension vector state by
    \begin{equation}
        \label{eqn:lorenz}
        \frac{du_j}{dt} = (u_{j+1} - u_{j-2})u_{j-1} - u_j + F, \quad {j = 1, \dots, 40}
    \end{equation}
Here, the variables are indexed periodically, e.g., $u_{41} = u_1$. The external forcing $F$ is the same for each mode, and so the equilibrium statistics of the system are spatially invariant. Note that quadratic nonlinearity satisfies the energy conservation law given in equation \eqref{eqn:energy_conservation}, which can be shown by the symmetry of the nonlinearity in equation \eqref{eqn:lorenz}.

A key property of the Lorenz '96 model is that it exhibits a variety of dynamical and statistical regimes by altering the value of $F$. Multiple dynamical regimes are typical of complex turbulent systems and represent a classic obstacle to effective control. Figure \ref{fig:lorenz_dynamics} exhibits two dynamical regimes corresponding to $F=5$ and $F=8$. The $F=5$ regime is weakly chaotic. It has highly non-Gaussian statistics and a long decorrelation time. Meanwhile, the regime corresponding to $F=8$ features near-Gaussian statistics and strongly chaotic dynamics with a correspondingly short decorrelation time. Previous results \cite{majda2017effective} have shown the statistical control strategy to be effective at controlling small perturbations in the Lorenz '96 model back to the equilibrium state.

\begin{figure}
        \begin{center}
            \includegraphics{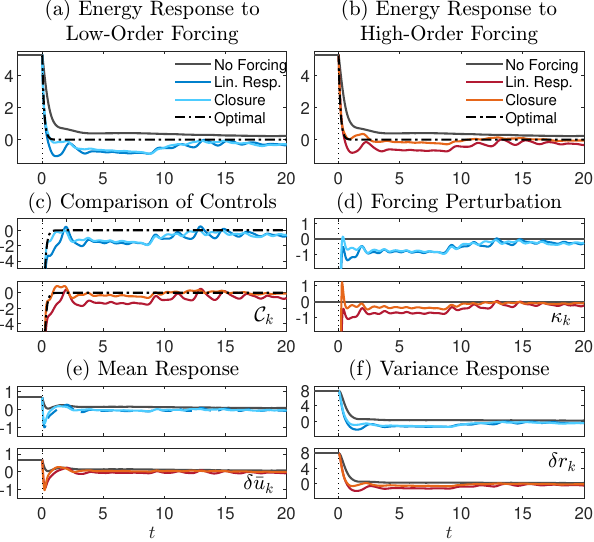}
        \end{center}
        \caption{\label{fig:lorenz}
            The control of the Lorenz 96 model from the perturbed state of $F = 8$ back to the equilibrium state of $F = 5$.
            Note this is a large perturbation into a regime with very different dynamics and statistics from the equilibrium.
            Panels (a) and (b) show the response of the energy perturbation to the control forcing for the low-order strategies and high-order strategies respectively.
            The energy perturbation is normalized by the dimension of the system.
            Panels (c)-(f) show the controls, forcing, mean response, and variance response for each mode.
            Note the system is translationally invariant, so the corresponding values for each mode are the same.
        }
    \end{figure}

In the current experiment, we show the efficacy of the strategies on a large perturbation which drive the system into a different dynamical regime. This leads to a much more challenging problem since the model state goes through a statistical transition between two distinctive regimes. The linear response estimation is no longer valid since the model moves far beyond the linear and near-Gaussian regime. The higher-order corrections become necessary to guarantee effective control performance. In this case, the $F_\mathrm{eq} = F = 5$ regime is taken as the equilibrium state and $\delta F_\mathrm{p}=3$ so that the perturbed state is in the $F=8$ regime. Figure \ref{fig:lorenz} shows the control strategies for this large perturbation. In this case, the high-order method with the mean dynamical equation closure shows the most effective strategy. In fact, it shows that combining both the high-order and the mean closure methods is essential to achieve good performance. This is a typical example to confirm the necessity of including high-order corrections when nonlinear and non-Gaussian features become dominant.

    \section{Further Discussions} \label{sec:discussion}

The statistical control strategies extend the effective control methods beyond the small perturbation scenario and demonstrate promise for applications to a broader range of turbulent situations. To summarize, they enjoy several attractive features. Using the total energy as the object of control, there is no need to track and control a large dimension of instabilities due to the energy-conservation principle. Thus, the computational cost is significantly reduced and independent of the dimensionality of the system. Further, the control can be determined entirely offline and only requires statistical information about the target equilibrium state, which is usually available from history observation data in many realistic applications.

Here, we discuss several key features in the new high-order control strategies based on the observations from the numerical experiments.

    \subsection{The High-Order Correction}

The control-forcing relation, which encodes the energy response to the external deterministic forcing perturbation, has a second-order perturbation term, $\delta \mathbf{F} \cdot \delta \bar{\mathbf{u}}$, which accounts for the higher-order contributions of the deterministic forcing perturbation to the energy response. This term consists of the product of the forcing perturbation and the mean perturbation in response to the forcing. Under the circumstances with small perturbations from the equilibrium state, both the mean perturbation and the forcing perturbation are small, so this term can be truncated without compromising the accuracy of the energy response. This method, where only the leading-order contributions to the energy response are considered, is the low-order method. However, for most large perturbations, the second-order term becomes large and significantly affects the energy response. The high-order approach incorporates this second-order term in the inversion of the control-forcing relation, fully resolving the energy response given the forcing perturbation and mean response. We explore several circumstances where the second-order perturbation term significantly impacts the energy response; thus, adding the high-order method can yield significant improvements over the low-order method.

When the initial mean perturbation is large, the initial value of the external forcing perturbation is greatly affected by the presence of the second-order term. This can be seen in equations \eqref{eqn:low_order_method} and \eqref{eqn:high_order_method} where the initial condition for the high-order method includes an extra term for the initial mean perturbation $\delta \bar{\mathbf{u}}(0)$. The effect of the initial mean perturbation on the resulting forcing perturbation can also be seen in $\kappa_2$ of Figure~\ref{fig:ep_regime} in the control of Regime I in Section \ref{sec:triad}. However, a large initial mean perturbation is unnecessary, and the second-order term can still have a significant effect even when the initial mean perturbation is relatively small if there is still a large initial energy perturbation due to the variance. Because the initial energy perturbation is large, the relative balance of the mean and variance in the total energy perturbation can shift over time, resulting in a potentially large mean perturbation after the initial time. In this case, the second-order term significantly impacts the evolution of the forcing perturbation even with the same initial conditions. Several examples of this phenomenon can be seen in the numerical tests in Sections \ref{sec:triad} and \ref{sec:lorenz}, especially in Figure~\ref{fig:lorenz}, where a drastic phase transition is shown.

As a further comment, whether to use the low-order or high-order methods can be made independently for each mode. For example, in a multiscale system, the effect of the second-order term may be small relative to the total energy response for small-scale modes and only be significant for larger-scale modes. In this case, the high-order method could be applied to only a subset of large-scale modes, while the low-order method is used for the rest, simplifying the dynamics in those small-scale modes without compromising performance.

    \subsection{The Mean Closure Equation}

The response of statistical energy to the external forcing depends directly on the mean response to the forcing. It is indirectly linked to the higher-order moments through the mean dynamical equation. This property is critical to formulating the statistical control strategy, allowing for attributing an external forcing perturbation to the optimal control by solving the control-forcing relation. Therefore, accurately approximating the mean response to external forcing is vital to the success of the statistical control strategy. Linear response theory effectively approximates the mean response for small perturbations, and it can perform well for larger perturbations in some cases when non-Gaussian statistics is not so important. However, its skill degenerates when the system is largely perturbed to a different dynamical regime where the linear response operator, based solely on the unperturbed dynamics, can provide very little information for the future perturbed state.

Using a mean closure equation for the mean response, which directly incorporates dynamics from the model, is expected to show improved performance when the initial perturbation spans multiple dynamical regimes. A mean dynamical closure equation is based on the explicit mean dynamics given in equation \eqref{eqn:meanDE}, in which the dependence on higher-order moments is closed using a suitable approximation. The closure considered in this paper utilizes a linear response for the higher-order contribution of the covariance described in equations \eqref{eqn:mean_closure} and \eqref{eqn:covariance_linear_response}. While this mean closure equation still relies on the linear response for the covariance, the mean dynamics still provide more information about the perturbed regime than the mean linear response.

The initial forcing perturbation is affected by the choice of mean response. The mean linear response cannot directly use the initial mean perturbation and instead must use the initial mean perturbation predicted by the linear response to a constant forcing perturbation. This is necessary to guarantee the convergence of the forcing perturbation to zero in the linear response case. The mean closure model, however, can utilize the initial mean perturbation directly. This is illustrated by $\kappa_1$ in Figure~\ref{fig:ep_regime} in Regime I of Section \ref{sec:triad} where the initial mean response differs between the linear response and mean closure methods. In this case, the mean equation closure achieves better performance. Even when the initial mean perturbation is similar to the initial mean perturbation predicted by the linear response, the mean closure model can provide more accurate dynamics in many cases. In Section \ref{sec:triad}, the mean equation closure method performs better among all cases, especially in Regime II. In Section \ref{sec:lorenz}, the Lorenz '96 model is perturbed from a non-Gaussian regime to a near-Gaussian turbulent regime. The mean equation closure again performs better than the linear response in this case with multiple dynamical regimes.

    \subsection{Convergence to the Equilibrium State}

\begin{figure}
        \begin{center}
        \includegraphics{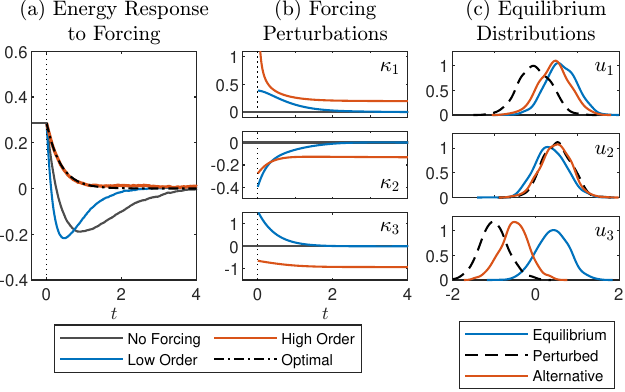}
        \end{center}

        \caption{\label{fig:opposite_regime}
            An example where the optimal energy response is achieved, but the system is forced to a different equilibrium state.
            The triad model has parameters $d_1=d_2=d_3=1$, $L_1=L_2=L_3=0$, $B_1=1$, $B_2=-0.6$, $B_3=-0.4$, $F_1=F_2=F_3=0.5$, and $\sigma_1=\sigma_2=\sigma_3=0.5$.
            The perturbed state has $F_3 = -1$.
            Panel (a) shows the energy response to the low-order and high-order control strategies.
            Panel (b) shows the forcing perturbations in each mode.
            Note that the forcing perturbation for the high-order method does not converge to zero.
            Panel (c) shows the equilibrium distribution, the perturbed distribution, and the alternative distribution achieved by the high-order method that yields the same statistical energy.
        }
    \end{figure}

The total statistical energy bounds the total mean and variance. Ideally, one hopes that the efficient control of the equilibrium energy will also achieve the efficient control of the mean and variance back to the equilibrium state. While this appears to be the case in most applications, this is not mathematically guaranteed by the statistical control strategy. Figure \ref{fig:opposite_regime} illustrates an example where the optimal energy response is achieved through the high-order mean closure strategy, but the external forcing perturbation converges to a constant non-zero state. Essentially the system converges to a different equilibrium with a different constant external forcing but the same equilibrium energy.

The cost functional given in equation \eqref{eqn:cost} only penalizes the strength of the direct energy control $\mathcal{C}_k$ rather than the external forcing perturbations $\kappa_k$ that yield that control. In addition, the control-forcing relation given in Equation \eqref{eqn:control_forcing_relation} admits multiple solutions in the limit in both the low-order and high-order formulations. One natural fix for such an issue is incorporating additional terms into the cost function, for example, explicitly excluding the mean state.

    \section{Conclusion} \label{sec:conclusions}

An efficient method of controlling the complex turbulent system with energy conserving nonlinearity is achieved through control of the total statistical energy  from a perturbed state back to equilibrium without controlling the large number of multiscale and potentially unstable modes. 
This paper proposed new statistical control strategies by extending previous works \cite{majda2017effective, majda2019using}, which had been restricted to scenarios with small perturbations from the equilibrium state. Incorporating the high-order term in the control-forcing relation accounts for the second-order contribution of the perturbations to the statistical energy, allowing the strategy to account for the response of the energy more accurately to large amplitude external forcing perturbations. Additionally, introducing a mean dynamical closure model allows the statistical control strategy to better account for large perturbations that drive the system into different dynamical regimes whose dynamics cannot be adequately reflected directly by a mean linear response approximation. These strategies allow for the practical application of the statistical control strategy to a wider variety of perturbations and regimes than previously possible.

The field of statistical control theory remains relatively underexplored, providing many promising research directions. The results presented in this paper could be further refined by developing more sophisticated methods for incorporating the mean and covariance dynamics into the mean response. Besides, other designs for mean closure models could be considered, several of which are described in \cite{majda2018strategies}. It may also be possible to incorporate other more suitable statistical functionals into the control strategy in addition to the statistical energy. This would allow, for example, the direct control of the statistical mean, giving finer control over the response of the system. Lastly, while the current statistical control strategy is conducted in an open-loop fashion, determining the control offline, a natural extension would be incorporating feedback from the system into a closed-loop statistical control strategy. This could be accomplished, for example, by combining measurements of the actual mean response of the system into the inversion of the control-forcing relation. The computational advantages of the statistical control strategy would be very advantageous in such a closed-loop formulation, which requires real-time incorporation of the model feedback.

    \section*{Acknowledgements}
The research of N.C. is funded by ONR N00014-21-1-2904 and ARO W911NF-23-1-0118. J.C. is partially supported by ONR N00014-21-1-2904 as a graduate research assistant. 
    \section*{Authors' Contributions}
J.C.: investigation, analysis, validation, visualization, writing—original draft; \\
D.Q.: methodology, supervision, analysis, writing—editing;\\
N.C.: methodology, supervision, analysis, writing—editing;\\
All authors gave final approval for publication and agreed to be held accountable for the work performed therein.
\section*{Conflict of interest declaration}
The authors have no conflict of interest.

    \nocite{gritsun2008climate}
    \nocite{bucci2019control}

    \nocite{majda2015statistical}

    \bibliography{refs}{}
    \bibliographystyle{plain}
\end{document}